\documentclass[12pt,oneside,leqno]{article}
\usepackage[T1]{fontenc}
\usepackage[latin1]{inputenc}
\usepackage{amsmath,latexsym}
\usepackage{amsfonts}
\newtheorem{theorem}{Theorem}
\newtheorem{lemma}{Lemma}
\newtheorem{proposition}{Proposition}
\newtheorem{definition}{Definition}

\begin{document}
\title{\normalsize\bfseries On the Irreducibility of the Complex Specialization of the Representation of
The Hecke Algebra of the Complex Reflection Group $G_7$ }
\author{\footnotesize MOHAMMAD Y. CHREIF AND MOHAMMAD N. ABDULRAHIM}
\date{}
\maketitle

\renewcommand{\thefootnote}{}

\noindent{abstract}
We consider a 2-dimensional representation of the Hecke algebra $\mathcal{H}(G_7, u)$, where $G_7$ is the complex reflection group
and $u$ is the set of indeterminates $u=(x_1, x_2, y_1, y_2, y_3, z_1, z_2, z_3)$. After specializing the indetrminates to non zero
complex numbers, we then determine a necessary and sufficient condition that guarantees the irreducibility of the complex specialization of the representation of 
the Hecke algebra $\mathcal{H}(G_7, u)$.

\begin{center}
\noindent {\footnotesize 2. INTRODUCTION }
\end{center}

Let $V$ be a complex vector space and $W$ a finite irreducible subgroup of $GL(V)$ generated by complex reflections. Let
$R$ be the set of reflections in $W$. For any element $s$ of $R$, denote by $H_s$ its pointwise fixed hyperplane. We define the set $V^{reg}=V$- $\cup_{s\in R}H_s$ and denote by $\bar{V}$ the quotient $V^{reg}/W$.

The braid group associated to $(W,V)$ is the fundamental group
$B(W)=\pi_1(\bar{V}, \bar{x}_0)$ of $\bar{V}$ with respect to any point $\bar{x}_0\in\bar{V}$.

We choose the set of indeterminates, $u=(u_{s,j})_{s,0\leq j\leq o(s)-1}$, where s runs over the generators of $W$ and $u_{s,j}=u_{t,j}$ if $s$ and $t$ are conjugate in $W$. Here $o(s)$ denotes the order of $s$. The cyclotomic Hecke algebra associated to $W$ is the quotient of the group algebra $\mathbb{Z}[u,u^{-1}]BW$ by the ideal generated by the relations
$\prod_{j=0}^{o(s)-1} (s-u_{s,j})$.

In [7], G. Malle and J. Michel constructed on the cyclotomic hecke algebra $\mathcal{H}(G_7, u)$ of the complex reflexion group, $G_7$, an 
irreducible representation $\phi:$ $\mathcal{H}(G_7, u)\rightarrow M_{2}(\mathbb{C}(u^{\frac{1}{2}},u^{-\frac{1}{2}}))$,
where $u$ is the set of indeterminates $u=(x_1, x_2, y_1, y_2, y_3, z_1, z_2, z_3)$.

In our work, we specialize the indeterminates $x_1, x_2, y_1, y_2, y_3, z_1$, $z_2$ and $z_3$ to nonzero complex numbers $\rho e^{i\alpha}$, where
$\alpha\in(-\pi, \pi]$ and $\rho$ a positive real number. We then get a representation $\varphi:$ $\mathcal{H}(G_7, u)\rightarrow GL_2(\mathbb{C})$. In section 3, we consider the case when $x_1=x_2$ and we show that $\varphi$ is irreducible if and only if
$z_1\neq\frac{y_1z_2}{y_2}$ and $z_1\neq\frac{y_2z_2}{y_1}$ (Theorem 2). In section 4, we assume that $x_1\neq x_2$ and we show that $\varphi$ is irreducible if and only if $x_1y_2z_2\neq x_2y_1z_1$,
$x_1y_1z_2\neq x_2y_2z_1$, $x_1y_2z_1\neq x_2y_1z_2$ and $x_1y_1z_1\neq x_2y_2z_2$ (Theorem 3).

\begin{center}
\noindent {\footnotesize 2. PRELIMINARIES}
\end{center}

\begin{definition}
\cite{G} Let $V$ be a complex vector space of dimension $n$. A complex reflection of GL(V) is a non-trivial element of GL(V) which acts trivialy on a hyperplane.
\end{definition}
\begin{definition}
\cite{G} Let $V$ be a complex vector space of dimension $n$. A complex reflection group is the subgroup of GL(V) generated by complex reflections.
\end{definition}

Examples of complex reflection groups include dihedral groups and symmetric groups. For $n\geq 3$, the dihedral group, $D_n$, is the group of the isometries of the plane preserving a regular
polygon, with the operation being composition.\\\\
A classification of all irreducible reflection groups shows that there are 34 primitive irreducible reflection groups [8].
The starting point was with A. Cohen, who provided a data for those irreducible complex reflection groups of rank 2 [5].

\begin{definition}
\cite{BMR} The complex reflection group, $G_7$, is an abstract group defined by the presentation
 \begin{center}
 $G_7=<t, u, s$ | $t^2=u^3=s^3=1, tus=ust=stu>$.
 \end{center}
\end{definition}

\begin{theorem}
\cite{B} The braid group of $G_7$ is isomorphic to the group

\begin{center}
$B=<s_1, s_2, s_3$ | $s_1s_2s_3=s_2s_3s_1=s_3s_1s_2>$.
\end{center}

\end{theorem}
Definitions and properties of braid groups are found in [2].
\begin{definition}
\cite{MG} Let $u$ be the set of indeterminates $u=(x_1, x_2, y_1, y_2, y_3, z_1, z_2, z_3)$. The cyclotomic Hecke algebra $\mathcal{H}(G_7, u)$ of $G_7$ is the quotient of the group algebra
of $B$ over $\mathbb{Z}[u, u^{-1}]$ by the relations

\begin{center}
$(s_1-x_1)(s_1-x_2)=0$,  \hspace{1cm}$\prod_{i=1}^{3}(s_2-y_i)=0$,  \hspace{1cm}$\prod_{i=1}^{3}(s_3-z_i)=0$. 
\end{center}
\end{definition}
For more details about the Hecke algebra of $G_7$, see [4].
\medskip
\begin{definition}
\cite{MG} Let $u=(x_1, x_2, y_1, y_2, y_3, z_1, z_2, z_3)$. The representation $\phi$ is defined as follows:
 
 \begin{center}
 $\phi: \mathcal{H}(G_7, u)\rightarrow M_{2}(\mathbb{C}(u^{\pm\frac{1}{2}}))$
 \end{center}
\vskip 1in

\begin{center}
$s_1=\begin{pmatrix}
x_1&\frac{y_1+y_2}{y_1y_2}-\frac{(z_1+z_2)x_2}{r}\\
0&x_2
\end{pmatrix},
\quad
s_2=\begin{pmatrix}
y_1+y_2&\frac{1}{x_1}\\
-y_1y_2x_1&0
\end{pmatrix}
$
\end{center}
and

\begin{center}
$s_3=\begin{pmatrix}
0&\frac{-r}{y_1y_2x_1x_2}\\
r&z_1+z_2
\end{pmatrix}$,
\end{center}
where $r=\sqrt{x_1x_2y_1y_2z_1z_2}$.
\end{definition}

We specialize the indeterminates $x_1, x_2, y_1, y_2, z_1, z_2$ and $z_3$ to nonzero complex numbers, $\rho e^{i\alpha}$, where
$\alpha\in(-\pi, \pi]$ and $\rho$ a positive real number. We then get a representation $\varphi:$ $\mathcal{H}(G_7, u)\rightarrow GL_2(\mathbb{C})$.
\begin{definition}
Principal square root function is defined as follows:\\
$z\in\mathbb{C}$, $z=(\rho,\alpha)$, $\rho\geq 0$. $\sqrt{z}=\sqrt{\rho}e^{i\frac{\alpha}{2}}$ where $-\pi<\alpha\leq\pi$.
\end{definition}
Since $\alpha\in(-\pi, \pi]$, it follows that $\sqrt{z^2}=z$ for any complex number z.
\begin{center}
\noindent {\footnotesize 3. IRREDUCIBILITY OF THE REPRESENTATION $\varphi$ FOR $x_1=x_2$ }
\end{center}

We assume that $x_1=x_2$ and we find a necessary and sufficient condition that guarantees the irreducibility of the representation
$\varphi$: $\mathcal{H}(G_7, u)\rightarrow GL_2(\mathbb{C})$.\\

Under this assumption, we have that the images of the generators of $\mathcal{H}(G_7, u)$ are\\
\begin{center}
$s_1=\begin{pmatrix}
x_2&\frac{y_1+y_2}{y_1y_2}-\frac{(z_1+z_2)x_2}{\sqrt{x_2^2y_1y_2z_1z_2}}\\
0&x_2
\end{pmatrix},\quad
s_2=\begin{pmatrix}
y_1+y_2&\frac{1}{x_2}\\-y_1y_2x_2&0
\end{pmatrix}
$
\end{center}
and
\begin{center}
$s_3=\begin{pmatrix}
0&-\frac{\sqrt{x_2^2y_1y_2z_1z_2}}{x_2^2y_1y_2}\\r&z_1+z_2
\end{pmatrix}.
$
\end{center}

For the matrix $s_1$, we denote by $s_1(i,j)$ the term of the matrix $s_1$ which lies in the $i$th row and in the $j$th column.\\

\begin{lemma}
$s_1(1,2)=0$ if and only if $z_1=\frac{y_1z_2}{y_2}$ or $z_1=\frac{y_2z_2}{y_1}$.
\end{lemma}
\textit{Proof.} We show that if $s_1(1,2)=0$ then $z_1=\frac{y_1z_2}{y_2}$ or $z_1=\frac{y_2z_2}{y_1}$.\\\\
Assume that $s_1(1,2)=0$. This implies that $\frac{y_1+y_2}{y_1y_2}=\frac{(z_1+z_2)x_2}{x_2\sqrt{y_1y_2z_1z_2}}$.
This implies that $(y_1+y_2)\sqrt{y_1y_2z_1z_2}=(z_1+z_2)y_1y_2$. Using $y_1y_2=(\sqrt{y_1y_2})^2$, we get $(y_1+y_2)\sqrt{z_1z_2}=(z_1+z_2)\sqrt{y_1y_2}$.\\\\
Squaring both sides, we obtain $(y_1+y_2)^2z_1z_2=(z_1+z_2)^2y_1y_2$. This implies that $z_1=\frac{y_1z_2}{y_2}$ or
$z_1=\frac{y_2z_2}{y_1}$.
On the other hand, direct computations show that if $z_1=\frac{y_1z_2}{y_2}$ or $z_1=\frac{y_2z_2}{y_1}$ then $s_1(1,2)=0$.

\vskip .5in

We now determine a sufficient condition for irreduciblity.
\begin{proposition}
The representation $\varphi$ is irreducible if $z_1\neq\frac{y_1z_2}{y_2}$ and $z_1\neq\frac{y_2z_2}{y_1}$.
\end{proposition}
\textit{Proof.}
Using the hypothesis and Lemma 1, we get $s_1(1,2)\neq 0$. Let $S$ be a non trivial proper invariant subspace of $\mathbb{C}^2$. The eigenspace of $s_1$ is generated by $e_1$. This implies that $S$ is of the form $<v>$,
where $v=ae_1$ for some non-zero complex number $a$.\\\\
$S$ is invariant implies that $s_2v=(a(y_1+y_2), -ax_2y_1y_2)\in S$, which is a contradiction. Therefore $S$ is irreducible.

\vskip .3in

We determine a necessary condition for irreducibility.
\begin{proposition}
The representation $\varphi$ is reducible if $z_1=\frac{y_1z_2}{y_2}$ or $z_1=\frac{y_2z_2}{y_1}$. 
\end{proposition}
\textit{Proof.} In each case, we show that the 1-dimensional subspace $M$ generated by the vector $u=(-\frac{1}{x_2y_2}, 1)$ is invariant.\\\\
\textbf{Case1.} $z_1=\frac{y_1z_2}{y_2}$. Substituting in Definition 5, we get

\begin{center}
$s_1=\begin{pmatrix}
x_2&0\\0&x_2
\end{pmatrix}, 
\quad\quad\quad\quad\quad\quad
s_2=\begin{pmatrix}
y_1+y_2&\frac{1}{x_2}\\-x_2y_1y_2&0
\end{pmatrix}
$
\end{center}
and
\begin{center}
$s_3=\begin{pmatrix}
0&-\frac{z_2}{x_2y_2}\\
x_2y_1z_2&z_2+\frac{y_1z_2}{y_2}
\end{pmatrix}.
$
\end{center}
\vskip .3in

It is easy to see that $s_2u=y_1u$ and $s_3u=z_2u$. This implies that $M$ is invariant.\\\\
\textbf{Case2.}  $z_1=\frac{y_2z_2}{y_1}$. Substituting in Definition 5, we get

\begin{center}
$s_1=\begin{pmatrix}
x_2&0\\0&x_2
\end{pmatrix},
\quad\quad\quad\quad\quad\quad
s_2=\begin{pmatrix}
y_1+y_2&\frac{1}{x_2}\\
-x_2y_1y_2&0
\end{pmatrix}
$
\end{center}
and

\begin{center}
$s_3=\begin{pmatrix}
0&-\frac{z_2}{x_2y_1}\\x_2y_2z_2&z_2+\frac{y_2z_2}{y_1}
\end{pmatrix}.
$
\end{center}
\vskip .2in

It is also easy to see that $s_2u=y_1u$ and $s_3u=z_2\frac{y_2}{y_1}u$. This implies that $M$ is invariant. 

\vskip .2in

Here we have proved the following theorem: 

\begin{theorem}
The representation $\varphi$ is irreducible if and only if $z_1\neq\frac{y_1z_2}{y_2}$ and $z_1\neq\frac{y_2z_2}{y_1}$.
\end{theorem}

\vskip .2in

\begin{center}
\noindent {\footnotesize 4. IRREDUCIBILITY OF THE REPRESENTATION $\varphi$ FOR $x_1\neq x_2$ }
\end{center}

We assume that $x_1\neq x_2$ and we find a necessary and sufficient condition that guarantees the irreducibility of the representation $\varphi$: $\mathcal{H}(G_7, u)\rightarrow GL_2(\mathbb{C})$.\\\\
For simplicity, we denote by $w$ the term\\\\
$(x_1-x_2)^2y_1^2y_2^2z_1z_2+[(y_1+y_2)r-x_1y_1y_2(z_1+z_2)][(y_1+y_2)r-x_2y_1y_2(z_1+z_2)]$...(1) 
\ \\
\begin{lemma}
The complex number $w$, defined in (1), is different from zero if and only if $x_1y_2z_2\neq x_2y_1z_1$, $x_1y_1z_2\neq x_2y_2z_1$, $x_1y_2z_1\neq x_2y_1z_2$ and $x_1y_1z_1\neq x_2y_2z_2$.
\end{lemma}
\textit{Proof.} Simple calculations show that $w=\alpha\beta$, where\\\\
$\alpha=x_2y_1y_2z_1+x_1y_1y_2z_2-(y_1+y_2)r$ and $\beta=x_1y_1y_2z_1+x_2y_1y_2z_2-(y_1+y_2)r$.\\\\
Assume that $w=0$. This implies that $\alpha=0$ or $\beta=0$.\\\\
If $\alpha=0$, then $x_2y_1y_2z_1+x_1y_1y_2z_2=(y_1+y_2)r$. Squaring both sides, we get
$y_1y_2(-x_2y_2z_1+x_1y_1z_2)(-x_2y_1z_1+x_1y_2z_2)=0$. This implies that $x_1y_1z_2=x_2y_2z_1$ or $x_1y_2z_2=x_2y_1z_1$.\\\\
If $\beta=0$, then $x_1y_1y_2z_1+x_2y_1y_2z_2=(y_1+y_2)r$. Squaring both sides, we get
$y_1y_2(x_1y_2z_1-x_2y_1z_2)(x_1y_1z_1-x_2y_2z_2)=0$. This implies that $x_1y_2z_1=x_2y_1z_2$ or $x_1y_1z_1=x_2y_2z_2$.\\\\
On the other hand, we assume that any of the following conditions holds true.
\begin{center}
$x_1y_2z_2=x_2y_1z_1$, $x_1y_1z_2=x_2y_2z_1$, $x_1y_2z_1=x_2y_1z_2$ or $x_1y_1z_1=x_2y_2z_2$
\end{center}
Under direct computations, we easily verify that $w=0$.

\vskip .2in

We now give a sufficient condition for the irreducibility of the representation $\varphi$.
\begin{proposition}
The representation $\phi$ is irreducible if $x_1y_2z_2\neq x_2y_1z_1$, $x_1y_1z_2\neq x_2y_2z_1$, $x_1y_2z_1\neq x_2y_1z_2$
and $x_1y_1z_1\neq x_2y_2z_2$.
\end{proposition}
\textit{Proof.} If the term $s_1(1,2)=\frac{y_1+y_2}{y_1y_2}-\frac{x_2(z_1+z_2)}{r}$ equals zero, then neither $e_1$ nor $e_2$ is a common eigenvector
for $s_2$ and $s_3$. This implies that the representation is irreducible. We note that under this case, we have that the complex number $w$ is not zero and hence, by Lemma 2, we also have that
$x_1y_2z_2\neq x_2y_1z_1$, $x_1y_1z_2\neq x_2y_2z_1$, $x_1y_2z_1\neq x_2y_1z_2$ and $x_1y_1z_1\neq x_2y_2z_2$.  \\\\
If $s_1(1,2)=\frac{y_1+y_2}{y_1y_2}-\frac{x_2(z_1+z_2)}{r}$ is not zero,
we diagonalize the matrix $S_1$ by the invertible matrix\\
\begin{center}
$T=\begin{pmatrix}
1&\frac{\frac{y_1+y_2}{y_1y_2}-\frac{x_2(z_1+z_2)}{r}}{x_2-x_1}\\
0&1
\end{pmatrix}.
$
\end{center}
We get 

\begin{center}
$T^{-1}s_1T=\begin{pmatrix}
x_1&0\\
0&x_2
\end{pmatrix}.
$
\end{center}
\ \\
We then conjugate $s_2$ by the matrix $T$. We get
\\
\begin{center}
$T^{-1}s_2T=\begin{pmatrix}
M&w\\
-x_1y_1y_2&P
\end{pmatrix},
$
where
\end{center}

\vskip .2in

$M=-\frac{x_2(-x_1y_1y_2z_1-x_1y_1y_2z_2+y_1r+ y_2r)}{(x_1-x_2)r}$,\\\\\\
and\\\\
$P=-\frac{x_1(-x_2y_1y_2z_1-x_2y_1y_2z_2-y_1r-y_2r)}{(x_1-x_2)r}$.\\\\

By conjugating $s_3$ by $T$, we get
\begin{center}
$T^{-1}s_3T=\begin{pmatrix}
A&B\\
r&C
\end{pmatrix},
$
where
\end{center}

\vskip .2in

$A=\frac{(y_1+y_2)r-x_2y_1y_2(z_1+z_2)}{(x_1-x_2)y_1y_2}$,\\\\
$B=\\\\\frac{1}{(x_1-x_2)^2r^3}x_1x_2z_1z_2(-x_1x_2y_1y_2z_1^2-x_1x_2y_1^2z_1z_2-x_1^2y_1y_2z_1z_2-2x_1x_2y_1y_2z_1z_2\\\\
-x_2^2y_1y_2z_1z_2-x_1x_2y_2^2z_1z_2-x_1x_2y_1y_2z_2^2+x_1y_1z_1r+x_2y_1z_1r\\\\+x_1y_2z_1r+x_2y_2z_1r+x_1y_1z_2r+x_2y_1z_2r
+x_1y_2z_2r+x_2y_2z_2r)$\\\\
and

\begin{center}
$C=\frac{-r(y_1+y_2)+x_1y_1y_2(z_1+z_2)}{(x_1-x_2)y_1y_2}$.
\end{center}

For simplicity, we denote $T^{-1}s_iT$ by $b_i$ for $1\leq i\leq 3$.\\\\

Suppose, to get contradiction, that the representation is reducible. That is, there exists a non trivial
proper invariant subspace $M$ of $\mathbb{C}^2$ of dimension 1.\\

The subspace $M$ has to be one of the following subspaces $<e_1>$ or $<e_2>$.\\\\
\textbf{Case1} $S=<e_1>$. Since $e_1\in M$, it follows that
 \begin{center}
 $b_3e_1=(A,r)\in M.$
 \end{center}
This implies that $r=0$, a contradiction.\\ \\
\textbf{Case2} $S=<e_2>$. Since $e_2\in M$, it follows that

\begin{center}
$b_2e_2=(w, P)\in M$.
\end{center}
By Lemma 2, we have that $w$ is different from zero, which is a contradiction.
Therefore the representation is irreducible.

We now present a lemma concerning the number $B$ used in defining $T^{-1}s_3T$ in Proposition 3.

\begin{lemma}
The complex number $B$ equals zero in each of the following cases:
\begin{enumerate}
\item $x_1y_2z_2=x_2y_1z_1$
\item $x_1y_1z_2=x_2y_2z_1$
\item $x_1y_2z_1=x_2y_1z_2$
\item $x_1y_1z_1=x_2y_2z_2$
\end{enumerate}
\end{lemma}
\textit{Proof.}  We verify that $B=0$ in case $(i)$. Suppose that $x_1y_2z_2=x_2y_1z_1$ then
$B=$\\\\
$\frac{1}{(x_1-x_2)^2r^3}x_1x_2z_1z_2(
-\frac{x_2y_1z_1}{y_2z_2}x_2y_1y_2z_1^2-\frac{x_2y_1z_1}{y_2z_2}x_2y_1^2z_1z_2-\frac{x_2^2y_1^2z_1^2}{y_2^2z_2^2}y_1y_2z_1z_2-2\frac{x_2y_1z_1}{y_2z_2}x_2y_1y_2z_1z_2\\\\
-x_2^2y_1y_2z_1z_2-\frac{x_2y_1z_1}{y_2z_2}x_2y_2^2z_1z_2-\frac{x_2y_1z_1}{y_2z_2}x_2y_1y_2z_2^2+
\frac{x_2y_1z_1}{y_2z_2}y_1z_1x_2y_1z_1+x_2y_1z_1x_2y_1z_1\\\\+\frac{x_2y_1z_1}{y_2z_2}y_2z_1x_2y_1z_1+x_2y_2z_1x_2y_1z_1+\frac{x_2y_1z_1}{y_2z_2}y_1z_2x_2y_1z_1+x_2y_1z_2x_2y_1z_1+\frac{x_2y_1z_1}{y_2z_2}y_2z_2x_2y_1z_1\\\\+x_2y_2z_2x_2y_1z_1$)\\\\
$=\frac{1}{(x_1-x_2)^2r^3}x_1x_2z_1z_2(-\frac{x_2^2y_1^2z_1^3}{z_2}-\frac{x_2^2y_1^3z_1^2}{y_2}-\frac{x_2^2y_1^3z_1^3}{y_2z_2}-2x_2^2y_1^2z_1^2-x_2^2y_1y_2z_1z_2-x_2^2y_1y_2z_1^2\\\\-x_2^2y_1^2z_1z_2+
\frac{x_2^2y_1^3z_1^3}{y_2z_2}+x_2^2y_1^2z_1^2+\frac{x_2^2y_1^2z_1^3}{z_2}+x_2^2y_1y_2z_1^2+\frac{x_2^2y_1^3z_1^2}{y_2}+x_2^2y_1^2z_1z_2+x_2^2y_1^2z_1^2\\\\+x_2^2y_1y_2z_1z_2)$\\\\
$=0$.\\\\
Likewise, we show that $B=0$ under each of the other conditions.

We now present a necessary condition for irreducibility.

\begin{proposition}
The representation is reducible in each of the following cases:
\begin{enumerate}
\item $x_1y_2z_2=x_2y_1z_1$
\item $x_1y_1z_2=x_2y_2z_1$
\item $x_1y_2z_1=x_2y_1z_2$
\item $x_1y_1z_1=x_2y_2z_2$
\end{enumerate}
\end{proposition}
\textit{Proof.} Assume that we have either one of the following conditions holds true:
\begin{center}
$x_1y_2z_2=x_2y_1z_1$, $x_1y_1z_2=x_2y_2z_1$, $x_1y_2z_1=x_2y_1z_2$ or $x_1y_1z_1=x_2y_2z_2$
\end{center}
Let $S$ be the one dimensional subspace generated by $e_2$.\\

If $s_1(1,2)=\frac{y_1+y_2}{y_1y_2}-\frac{x_2(z_1+z_2)}{r}$ equals zero, then $w$, as defined in section 4, equals
 $(x_1-x_2)^2y_1^2y_2^2z_1z_2$. This implies that $w\neq0$. By lemma 2, we get a contradiction.\\

Therefore, without loss of generality, we assume that $s_1(1,2)\neq0$. We then conjugate the representation by the invertible matrix $T$. Recall that $b_i=T^{-1}s_iT$ ($i=1, 2, 3$). We then
have that $b_2e_2=(w, P)=(0, P)$ by Lemma 2, and $b_3e_2=(B, C)=(0,C)$ by Lemma 3. It follows that $S$ is invariant under this representation.

This leads us to state a necessary and sufficient condition for the irreducibility of the representation.

\begin{theorem}
The representation is irreducible if and only if $x_1y_2z_2\neq x_2y_1z_1$, $x_1y_1z_2\neq x_2y_2z_1$, $x_1y_2z_1\neq x_2y_1z_2$
and $x_1y_1z_1\neq x_2y_2z_2$.
\end{theorem}

{\bf Conflict of Interest}

The authors declare that there is no conflict of interest regarding the publication of this paper.

\textit{Mohammad Y. Chreif, Department of Mathematics, Beirut Arab University, P.O. Box 11-5020, Beirut, Lebanon\\
E-mail address: myc102@student.bau.edu.lb}\\

\textit{Mohammad N. Abdulrahim, Department of Mathematics, Beirut Arab University, P.O. Box 11-5020, Beirut, Lebanon\\
E-mail address: mna@bau.edu.lb}


\vskip .2in

\begin{thebibliography}{1}
\bibitem{B} D. Bessis, J. Michel, {\em Explicit presentations for exceptional braid groups}, Experiment. Math., 13 (3), 
257-266, 2004.
\bibitem{Bir} J. Birman, {\em Braids, Links and Mapping Class Groups}, Annals of Mathematical Studies, Princeton University Press, 
82, 1975. 
\bibitem{BMR} M. Brou$\acute{e}$, G. Malle, R. Rouquier, {\em Complex reflection groups, braid groups, Hecke algebras}, J.
reine angew. Math., 500, 127-190, 1998.
\bibitem{Ch} M. Chlouveraki, {\em Degree and Valuation of the Schur elements of cyclotomic Hecke algebras}, J. Algebra, 320 (11),
3935-3949, 2008.
\bibitem{C} A. Cohen, {\em Finite complex reflection groups}, Ann. Sci. $\acute{E}$cole Norm. Sup. (4), 9 (3), 379-436, 1976.
\bibitem{G} I. Gordon, S. Griffeth, {\em Catalan numbers for complex reflection groups}, Amer. J. Math., 134 (6), 1491-1502, 2012.
\bibitem{MG} G. Malle, J. Michel, {\em Constructing representations of Hecke algebras for complex reflection groups}, LMS J. Comput. Math., 13, 426-450, 2010.
\bibitem{S} G. Shephard, J. Todd, {\em Finite unitary reflection groups}, Canadian J. Math., 6, 274-304, 1954.
\end{thebibliography}
\end{document}